\documentclass[11pt]{article}
\usepackage{amsfonts}
\usepackage{amsfonts}
\usepackage{amsfonts}
\usepackage{amsfonts}
\usepackage{amsfonts}
\usepackage{amscd}
\usepackage{amsmath}
\usepackage[all]{xy}
\usepackage{latexsym}
\usepackage{amssymb}
\usepackage{amsfonts}
 \usepackage{mathrsfs}
 \usepackage[all]{xy}
 \usepackage{latexsym,amsfonts,amssymb,bm}

\begin{document}
\baselineskip=0.6truecm

\font\one=cmbx10 scaled\magstep4 \font\fif=cmti10 scaled\magstep1
\font\fiv=cmti10 scaled\magstep5 \font\ooo=cmsy10 scaled\magstep3
\font\two=cmcsc10 \font\three=cmti8 \font\ss=cmsl8 \font\bb=cmbx8
\font\four=cmbx10 scaled\magstep1 \font\six=cmti10 scaled\magstep1
\font\small=cmr8

\def\f{\noindent}
\def\c{\centerline}
\def\ol{\overline}
\def\T{\text{T}}
\def\t{\text}
\def\d{{\text{d}}}
\def\om{\omega}
\def\Om{\Omega}
\def\sub{\subset}
\def\al{\alpha}
\def\dt{\delta}
\def\ep{\varepsilon}
\def\eq{\equiv}
\def\la{\lambda}
\def\lg{\langle}
\def\rg{\rangle}
\def\e{{{\text{e}}}}
\def\R{{\Bbb R}}
\def\C{{\Bbb C}}
\def\A{{\cal A}}
\def\B{{\cal B}}
\def\vp{\varphi}
\def\pt{\partial}
\def\na{\nabla}
\def\Dt{\Delta}
\def\pma{\pmatrix}
\def\epm{\endpmatrix}
\def\Ga{\Gamma}
\def\for{\forall}
\def\si{\sigma}
\def\wt{\widetilde}
\def\tha{\theta}
\def\V{{\cal V}}
\def\be{\beta}
\def\bm{\bmatrix}
\def\ebm{\endbmatrix}
\def\ga{\gamma}
\def\us{\underset}
\def\os{\overset}
\def\1{{\pmb 1}}
\def\tri{\triangle}
\def\Si{\Sigma}
\def\U{{\cal U}}

\title{Pr\"{u}fer sheaves and generic sheaves
over the weighted projective lines and elliptic curves
\thanks{Supported in part by the National Natural
 Science Foundation of China (Grant No. 11201386), the Natural Science Foundation of Fujian Province of China (Grant No.
 2012J05009) and the Fundamental Research Funds for the Central Universities of China(Grant No.
 2012121004)}}
\author{\small\footnotesize  Jianmin Chen, Jinjing Chen, Yanan Lin
\thanks{  Corresponding author, $E-mail:$ jinjingchenyu@126.com}\\
  {\small\footnotesize The School of Mathematical
Sciences, Xiamen
 University, Xiamen 361005,  P.R. China}
 }

\date{}
\maketitle \par \noindent {\bf Abstract:} \ In the present paper, we
introduce the concepts of Pr\"{u}fer sheaves and adic sheaves over a
weighted projective line of genus one or an elliptic curve, show that
Pr\"{u}fer sheaves and adic sheaves can classify the category of coherent sheaves. Moreover, we describe the relationship
between Pr\"{u}fer sheaves and generic sheaves, and provide two methods to construct generic sheaves by
using coherent sheaves and Pr\"{u}fer sheaves.\\
  {\bf Keywords}: \  weighted projective lines; elliptic curves; quasi-coherent sheaves; pure-injective objects.\\
  {\bf 2000 Mathematics Subject Classification:} \ 14F05, 16G10, 16G60, 16G70, 18A30\\

\section{Introduction}

\hspace{4mm} The notion of weighted projective lines was introduced
by Geigle and Lenzing [10] to give a geometric treatment to canonical
algebras which was studied by Ringel [18]. Let $k$ be an
algebraically closed field, a weighted projective line over $k$ can
be viewed as obtained from a projective line $\mbox{P}_{k}^{1}$ by
endowing with positive integral multiplicities $p_{1},$ \ldots
$,p_{t}$( which were called weights) on pairwise distinct points
$\lambda_{1},$ \ldots $,\lambda_{t}$. An elliptic curve over $k$ is a smooth plane projective curve of genus one having a k-rational point. Much interesting work has been
done on the category of coherent sheaves over a weighted projective
line or an elliptic curve. But we still know little about the category of quasi-coherent
sheaves over them, even the properties of
some special quasi-coherent sheaves.

As we know that, quasi-coherent sheaves and coherent sheaves over a
scheme play the similar role with modules and finitely generated
modules over rings. And in the representation theory of algebras, some
special modules, for example generic modules, Pr\"{u}fer modules,
adic modules, have been widely studied ([3], [5], [17], [19], [20]).
These infinitely dimensional modules play an important role in the category of modules.
In 1997, Lenzing [14]
extended the concept of generic modules to generic sheaves over
weighted projective lines of genus one and determined all
indecomposable generic sheaves. Moreover, he proved that the left
perpendicular category of the generic sheaf $G_{q}$ with $q\in
\mathbb{Q}\bigcup\{\infty\}$ intersecting the category of coherent
sheaves is exactly the subcategory $\mathcal{C}^{(q)}$.

In view of the important role of generic modules, Pr\"{u}fer modules
and adic modules in the representation theory of finite-dimensional algebras
and generic sheaves over weighted projective lines of genus one, we
make further study about two special quasi-coherent sheaves which we
call Pr\"{u}fer sheaves and adic sheaves in this paper. The paper is
organized as follows:

In section 2, we recall the structure of the category of coherent sheaves over a weighted projective
line of genus one or an elliptic curve. In section 3, we extend the concepts of Pr\"{u}fer modules and adic
modules to Pr\"{u}fer sheaves and adic sheaves over a weighted
projective line of genus one or an elliptic curve and show how to use
Pr\"{u}fer sheaves and adic sheaves to classify the category of coherent sheaves. In section 4,
we prove that generic sheaves, Pr\"{u}fer sheaves and adic sheaves
are pure-injective
objects.
Section 5 describes an important relationship between Pr\"{u}fer sheaves, adic
sheaves and generic
sheaves in Theorem 5.3 and Theorem 5.4.
We provide two methods to construct
generic sheaves by using Pr\"{u}fer sheaves and coherent sheaves in
section 6.

In the paper, $\mbox{Hom}(X,Y)$ means
$\mbox{Hom}_{\mbox{Qcoh}(\mathbb{X})}(X,Y)$,
$\mbox{Ext}^{1}(X,Y)$ means $\mbox{Ext}^{1}_{\mbox{Qcoh}(\mathbb{X})}(X,$
$Y)$, for each $X,Y\in
\mbox{Qcoh}(\mathbb{X}).$ And we view the isomorphism as equality.

\section{The category of coherent sheaves on a weighted projective line of genus one or an elliptic curve}

Let $k$ be an algebraically closed field, $\mathbb{X}$
be a weighted projective line of genus one over $k$. It is
well-known that every quasi-coherent sheaf on $\mathbb{X}$ is a
direct limit of coherent sheaves, and the category
$\mbox{Qcoh}(\mathbb{X})$ of quasi-coherent sheaves on $\mathbb{X}$ is
a locally noetherian Grothendieck category. Hence, the structure of
a quasi-coherent sheaf on $\mathbb{X}$ much depend on that of
coherent sheaves on $\mathbb{X}$.  In this section, we recall the
structure of the category of coherent sheaves on $\mathbb{X}$.

{\bf{Proposition 2.1}} (see [10]) The category
$\mbox{coh}(\mathbb{X})$ of coherent sheaves on $\mathbb{X}$ is an
abelian, Ext-finite, noetherian, hereditary and Krull-Schmidt
$k$-category. $\mbox{coh}(\mathbb{X})$ satisfies Serre duality, i.e.
for any two coherent sheaves $F$ and $G$, there is an isomorphism
$\mbox{Hom}(F, \tau{G})= \mbox{DExt}^{1}(G, F)$, where $\mbox{D} =
\mbox{Hom}_{k}(-, k)$.

In addition, $\mbox{coh}(\mathbb{X})=
\mbox{coh}_{+}(\mathbb{X})\bigvee \mbox{coh}_{0}(\mathbb{X})$, that
is, each indecomposable object of $\mbox{coh}(\mathbb{X})$ lies
either in $\mbox{coh}_{+}(\mathbb{X})$ or in
$\mbox{coh}_{0}(\mathbb{X})$, and there are no non-zero morphisms
from $\mbox{coh}_{0}(\mathbb{X})$ to $\mbox{coh}_{+}(\mathbb{X})$,
where $\mbox{coh}_{+}(\mathbb{X})$ denotes the full subcategory of
$\mbox{coh}(\mathbb{X})$ consisting of all objects which do not have
a simple subobject, and $\mbox{coh}_{0}(\mathbb{X})$ denotes the
full subcategory of $\mbox{coh}(\mathbb{X})$ consisting of all
objects of finite length.

For more detail structure of $\mbox{coh}(\mathbb{X})$, we need
introduce rank, degree and slope of coherent sheaves.

Let $F, G\in \mbox{coh}(\mathbb{X})$, the Euler form of $F$ and $G$
is defined by
$$\langle F,G\rangle=\mbox{dim}_{k}\mbox{Hom}(F,
G)-\mbox{dim}_{k}\mbox{Ext}^{1}(F,G).$$ \noindent which can induce a
non-degenerated bilinear form $\langle-,-\rangle$ on the
Grothendieck group $\mathrm{K}_{0}(\mathbb{X})$, also called Euler
form. We call $x\in\mbox{K}_{0}(\mathbb{X})$ a radical vector if $x$ lies in
the radical $\mbox{K}_{0}(\mathbb{X})$ of the associated quadratic form.

{\bf{Lemma 2.2}} The radical of the Grothendieck group
$\mathrm{K}_{0}(\mathbb{X})$ has a $\mathbb{Z}$-basis $u$, $w$ such
that $\langle u, w \rangle= p =\mbox{l.c.m.}(p_{1},\ldots, p_{t})$, where
$(p_{1},\ldots, p_{t})$ is the weight sequence of $\mathbb{X}$.

{\bf{Definition 2.3}} For each coherent sheaf $F$, define the rank
of $F$ by $\mbox{rk}(F)=\langle [F], w\rangle$, and the
degree of $\mathcal{F}$ by $\mbox{deg}(F)=\langle u, [F]\rangle$,
where $[F]\in \mathrm{K}_{0}(\mathbb{X})$ is the corresponding class
of $F$. Then the slope of a coherent sheaf $F$ is an element in
$\mathbb{Q}\bigcup\{\infty\}$ defined as $\mu(F) =
\mbox{deg}(F)/\mbox{rk}(F)$.

{\bf{Proposition 2.4}} (see [15]) For each
$q\in\mathbb{Q}\bigcup\{\infty\}$, let $\mathcal{C}^{(q)}$ be the
additive closure of the full subcategory of $\mbox{coh}(\mathbb{X})$
formed by all indecomposable coherent sheaves of slope $q$. Then the
following holds:

(i) $\mathcal{C}^{(q)}$ is isomorphic to
$\mbox{coh}_{0}(\mathbb{X})$ for each
$q\in\mathbb{Q}\bigcup\{\infty\}$. In particular,
$\mathcal{C}^{(\infty)}$ is just $\mbox{coh}_{0}(\mathbb{X})$, which
is uniserial, i.e. each indecomposable object has a unique finite
composition series, and admits a natural decomposition
$\mbox{coh}_{0}(\mathbb{X})=\coprod_{x\in\mathbb{X}}U_{x}$, where
$U_{x}$ are connected uniserial categories indexed by $\mathbb{X}$.

(ii) $\mbox{coh}(\mathbb{X})$ is the additive closure of
$\bigcup_{q\in\mathbb{Q}\bigcup\{\infty\}}\mathcal{C}^{(q)}$.

(iii) $\mbox{Hom}_{\mathbb{X}}(\mathcal{C}^{(q)},
\mathcal{C}^{(r)})\neq 0$ if and only if $q\leq r$.

(iv) (Riemann-Roch formula) \ For each $F,G\in
\mbox{coh}(\mathbb{X})$, there has
$$\sum\limits_{i=0}^{p-1}\langle[\tau^{i}F],[G]\rangle=\mbox{rk}(F)\mbox{deg}
(G)-\mbox{rk}(G)\mbox{deg}(F),$$
where $ p =\mbox{l.c.m.}(p_{1},\ldots, p_{t})$.

An elliptic curve $\mathbb{E}$ over $k$ is a smooth plane projective curve of genus one having a k-rational point, the category $\mbox{coh}(\mathbb{E})$ of coherent sheaves has a similar structure as the weighted projective line of genus one where Proposition 2.1 and 2.4 also hold in $\mathbb{E}$. More detail of elliptic curves can be referred in [1], [2], [4], [13].

\section{Pr\"{u}fer sheaves and adic sheaves}

Let $\mathbb{X}$ be a weighted projective line of genus one or an elliptic curve. In this section, we introduce the concepts of Pr\"{u}fer sheaves and adic sheaves over $\mathbb{X}$, and discuss the set of morphisms between these two classes of quasi-coherent sheaves and coherent
sheaves.

By Proposition 2.4, we know that for each
$q\in\mathbb{Q}\bigcup\{\infty\}$, the Auslander-Reiten quiver of
$\mathcal{C}^{(q)}$ consists of stable tubes indexed by
$\mathbb{X}$.

{\bf{Definition 3.1}} Let $\mathcal{T}$ be a stable tube in
$\mathcal{C}^{(q)}$ with the rank $d$. Let $S_{q}$ be a quasi-simple
sheaf (i.e., simple object in $U_{x}$ whose Auslander-Reiten quiver is
$\mathcal{T}$) belonging to $\mathcal{T}$ and $S_{q}[i]$ be the
indecomposable sheaf of length $i$ in $\mathcal{C}^{(q)}$ satisfying
$\mbox{Hom}(S_{q}, S_{q}[i])\neq 0$. Then there is a sequence of embeddings
$$S_{q}\rightarrow S_{q}[2]\rightarrow\ldots\rightarrow
S_{q}[i]\rightarrow\ldots.$$ Denote by $S_{q}[\infty]$ the corresponding
direct limit. Composing the irreducible morphisms between the
sheaves belonging to $\mathcal{T}$ in the appropriate way we obtain
a generalized tube $\mathcal{T}'=(S_{q}[di+1])_{i\in\mathbb{N}_{0}}$.
Comparing to the definition of Pr\"{u}fer modules in [12],  we call
$S_{q}[\infty]$ a Pr\"{u}fer sheaf of slope $q$ over $\mathbb{X}$.

Similarly, there is also an indecomposable sheaf $S_{q}[-i]$ of length
$i$ in $\mathcal{C}^{(q)}$ satisfies $\mbox{Hom}(S_{q}[-i],$
$S_{q})\neq 0$, and we can obtain a sequence of epimorphisms
$$\ldots\rightarrow S_{q}[-i]\rightarrow\ldots\rightarrow S_{q}[-2]\rightarrow S_{q}.$$ Denote by
$S_{q}[-\infty]$ the corresponding inverse limit. We call $ S_{q}[-\infty]$
an adic sheaf of slope $q$ over $\mathbb{X}$.

Next we talk about
the sets of morphisms between coherent sheaves
and Pr\"{u}fer sheaves, and then between coherent sheaves and adic sheaves.
We need the following lemmas.

{\bf{Lemma 3.2}} (i) Let $X\in \mbox{coh}\mathbb{(X)}$. If $\{Y_{i}\ |\
i\in I, Y_{i}\in \mbox{Qcoh}\mathbb{(X)}\}$ is a direct system, then
\begin{center}
$\mbox{Hom}(X,
\underrightarrow{\mbox{lim}}Y_{i})=\underrightarrow{\mbox{lim}}\mbox{Hom}(X,
Y_{i})$  \ \ and   \ \ $\mbox{Ext}^{1}(X,
\underrightarrow{\mbox{lim}}Y_{i})=\underrightarrow{\mbox{lim}}\mbox{Ext}^{1}(X,
Y_{i})$.
\end{center}
If $\{Y_{i}\ |\ i\in I, Y_{i}\in \mbox{Qcoh}\mathbb{(X)}\}$ is a
inverse system, then
\begin{center}
$\mbox{Hom}(X,
\underleftarrow{\mbox{lim}}Y_{i})=\underleftarrow{\mbox{lim}}\mbox{Hom}(X,
Y_{i}).$
\end{center}

(ii) Let $X, Y\in\mbox{Qcoh}{(\mathbb{X})}$. If $Y=\underrightarrow{\mbox{lim}}Y_{i}$ with $Y_{i}\in\mbox{Qcoh}(\mathbb{X})$ and
$\mbox{Ext}^{1}(Y_{i}, X)$

\noindent $=0$ for every $i$, then $\mbox{Ext}^{1}(Y, X)=0$.
Dually, if $Y=\underleftarrow{\mbox{lim}}{Y_{i}}$ with $Y_{i}\in\mbox{Qcoh}({\mathbb{X}})$
and $\mbox{Ext}^{1}(X, Y_{i})=0$ for every $i$, then $\mbox{Ext}^{1}(X, Y)=0$.

$\bf{Proof:} $ (i) We only prove $\mbox{Ext}^{1}(X,
\underrightarrow{\mbox{lim}}Y_{i})=\underrightarrow{\mbox{lim}}\mbox{Ext}^{1}(X,
Y_{i})$, the rest formulas are obvious. For $\mbox{Qcoh}(\mathbb{X})$
has enough injective objects, we have an exact sequence
$$0\rightarrow Y_{i}\rightarrow I_{i}\rightarrow H_{i}\rightarrow 0,
\ \ \mbox{ for each}  \ Y_{i}$$ \noindent where $I_{i}$ is
injective. Applying $\mbox{Hom}(X, -)$ on these exact sequences, we
have long exact sequences
\begin{center}
$\mbox{Hom}(X, I_{i})\rightarrow \mbox{Hom}(X, H_{i})\rightarrow
\mbox{Ext}^{1}(X, Y_{i})\rightarrow 0,$
\end{center}
\noindent Taking direct limit, we have
\begin{equation}
\underrightarrow{\mbox{lim}}\mbox{Hom}(X, I_{i})\rightarrow
\underrightarrow{\mbox{lim}}\mbox{Hom}(X, H_{i})\rightarrow
\underrightarrow{\mbox{lim}}\mbox{Ext}^{1}(X, Y_{i})\rightarrow 0.
\end{equation}

\noindent On the other hand, applying direct limit on
$0\rightarrow Y_{i}\rightarrow I_{i}\rightarrow H_{i}\rightarrow 0$, we
obtain new exact sequences $0\rightarrow \underrightarrow{\mbox{lim}}Y_{i}
\rightarrow \underrightarrow{\mbox{lim}}I_{i}\rightarrow
\underrightarrow{\mbox{lim}}H_{i}\rightarrow 0$. Applying $\mbox{Hom}(X, -)$ to
it, we have long exact sequences
\begin{center}
Hom$(X, \underrightarrow{\mbox{lim}}I_{i})\rightarrow$ Hom$(X,
\underrightarrow{\mbox{lim}} H_{i}) \rightarrow$ Ext$^{1}(X,
\underrightarrow{\mbox{lim}} Y_{i})\rightarrow$ Ext$^{1}(X,
\underrightarrow{\mbox{lim}} I_{i}).$
\end{center}
$\underrightarrow{\mbox{lim}} I_{i}$ is injective since $\mbox{Qcoh}\mathbb{(X)}$
is hereditary, so we obtain another long exact sequence
\begin{equation}
\mbox{Hom}(X, \underrightarrow{\mbox{lim}}I_{i})\rightarrow
\mbox{Hom}(X, \underrightarrow{\mbox{lim}} H_{i}) \rightarrow
\mbox{Ext}^{1}(X, \underrightarrow{\mbox{lim}} Y_{i})\rightarrow 0.
\end{equation}
Compare (3.1) with (3.2), by the five lemma, we have $\mbox{Ext}^{1}(X,
\underrightarrow{\mbox{lim}}Y_{i})=\underrightarrow{\mbox{lim}}\mbox{Ext}^{1}(X,
Y_{i})$.

(ii) We prove the first result, the rest is duality. Let $0\rightarrow X\rightarrow Z\rightarrow Y\rightarrow 0$ be an exact sequence. Since $Y=
\underrightarrow{\mbox{lim}}Y_{i}$ and $\mbox{Ext}^{1}(Y_{i}, X)=0$, there exists an exact commutative diagram
\[\xymatrix{0\ar[r]&X\ar[r]\ar@{=}[d]&Z_{i}\ar^{\pi_{i}}[r]\ar^{\psi_{i}}[d]&Y_{i}\ar[r]\ar^{\phi_{i}}[d]&0\\
0\ar[r]&X\ar[r]\ar@{=}[d]&Z_{j}\ar^{\pi_{j}}[r]\ar[d]&Y_{j}\ar[r]\ar[d]&0\\
0\ar[r]&X\ar[r]&Z\ar^{\pi}[r]&Y\ar[r]&0}\]
for each $i<j$, where $\pi_{i}=(0, 1)$ and $\psi_{i}=
\left(
\begin{array}{cc}
 1 & 0\\
 0 & \phi_{i}
\end{array}
\right)
$.
Moreover, there exists a morphism $\eta_{i}=\left(
\begin{array}{c}
  0\\
  1
\end{array}
\right): Y_{i}\rightarrow Z_{i}$
satisfying $\pi_{i}\eta_{i}=1_{Y_{i}}$ and $\psi_{i}\eta_{i}=\eta_{j}\phi_{i}$. Then by the universality of direct limit, it induces a morphism $\eta:
Y\rightarrow Z$, we can prove that $\pi\eta=1_{Y}$. It implies $\mbox{Ext}^{1}(Y, X)=0$.
$\hfill\square$

{\bf{Lemma 3.3}} Let $X\in \mbox{coh}(\mathbb{X}), Y\in
\mbox{Qcoh}(\mathbb{X})$. There is an Auslander-Reiten formula
\begin{center}
 $\mbox{DExt}^{1}(X, Y)=\mbox{Hom}(Y, \tau{X})$.
\end{center}

$\bf{Proof:} $  Noticing that each quasi-coherent sheaf on
$\mathbb{X}$ is a direct limit of its coherent subsheaves, we may
assume that $Y=\underrightarrow{\mbox{lim}} Y_{i}$, where $\{Y_{i}\
|\ i\in I, Y_{i}\in \mbox{coh}\mathbb{(X)}\}$ is a direct system.
Using Lemma 3.2 and by Serre duality, we have $\mbox{DExt}^{1}(X,
Y)=\mbox{DExt}^{1}(X,
\underrightarrow{\mbox{lim}}Y_{i})=\mbox{D}\underrightarrow{\mbox{lim}}\mbox{Ext}^{1}(X,
Y_{i})=\underleftarrow{\mbox{lim}}\mbox{DExt}^{1}(X,
Y_{i})=\underleftarrow{\mbox{lim}}\mbox{Hom}(Y_{i}, \tau
X)=\mbox{Hom}(Y, \tau X)$.   $\hfill\square$

{\bf{Proposition 3.4}} Let $S_{q}[\infty]$ be the Pr\"{u}fer sheaf
of slope $q$ and $E$ be an indecomposable coherent sheaf which lies
in the mouth of a tube.

(i) If $\mu(E)<q$, then $\mbox{Ext}^{1}(E, S_{q}[\infty])=0$ and
$\mbox{Hom}(E, S_{q}[\infty])\neq 0$.

(ii) If $\mu(E)=q$, then $\mbox{Ext}^{1}(E,
S_{q}[\infty])=0=\mbox{Hom}(E, S_{q}[\infty])$ when $E\neq S_{q}$,
otherwise $\mbox{Ext}^{1}(E, S_{q}[\infty])=0$ and $\mbox{Hom}(E,
S_{q}[\infty])\neq 0$. In particular, $\mbox{dim}_{k}\mbox{Hom}(S_{q}, S_{q}[\infty])
=1$.

(iii) If $\mu(E)>q$, then $\mbox{Ext}^{1}(E, S_{q}[\infty])\neq 0$
and $\mbox{Hom}(E, S_{q}[\infty])=0$.

$\bf{Proof:}$
Assume that $S_{q}$ lies in a tube with rank $d$.

(i) If $\mu(E)<q$, we have
 $\mbox{Ext}^{1}(E,
S_{q}[\infty])=\underrightarrow{\mbox{lim}}\mbox{Ext}^{1}(E,
S_{q}[i])=0,$ and $\mbox{Hom}(E,$$S_{q}[\infty])=\underrightarrow{\mbox{lim}}\mbox{Hom}(E,
S_{q}[i])$. By Riemann-Roch formula, $\mbox{Hom}(E,$ 

\noindent $S_{q}[id])\neq 0$ for $i\in\mathbb{N}$. So $\mbox{Hom}(E, S_{q}[\infty])\neq0$.

(ii) If $\mu(E)=q$, the result is obvious by Auslander-Reiten formulas.

(iii) If $\mu(E)>q$, then $\mbox{Hom}(E, S_{q}[\infty])=0$. By
Riemann-Roch formula, we know that $\mbox{Hom}(S_{q}[id], \tau
E)\neq 0$ for $i\in \mathbb{N}$. Moreover, a non zero morphism from
$S_{q}[id]$ to $\tau E$ can be extended to a non zero morphism from
$S_{q}[\infty]$ to $\tau E$, so by Lemma
3.3, we have $\mbox{Ext}^{1}(E, S_{q}[\infty])\neq
0$.$\hfill\square$

{\bf{Proposition 3.5}} Let $S_{q}[-\infty]$ be the adic sheaf of
slope $q$ and $E$ be an indecomposable coherent sheaf lies in the
mouth of a tube.

(i) If $\mu(E)<q$, then $\mbox{Hom}(E, S_{q}[-\infty])\neq 0$ and
$\mbox{Ext}^{1}(E, S_{q}[-\infty])=0$.

(ii) If $\mu(E)=q$, then $\mbox{Hom}(E,
S_{q}[-\infty])=0=\mbox{Ext}^{1}(E, S_{q}[-\infty])$ when $E\neq
\tau^{-1}S_{q}$, otherwise $\mbox{Hom}(E, S_{q}[-\infty])=0$ and
$\mbox{Ext}^{1}(E, S_{q}[-\infty])\neq 0$. In particular, we have
$\mbox{dim}_{k}\mbox{Ext}^{1}(\tau^{-1}S_{q}, S_{q}[-\infty])$
$=1$.

(iii) If $\mu(E)>q$, then $\mbox{Hom}(E, S_{q}[-\infty])=0$ and
$\mbox{Ext}^{1}(E, S_{q}[-\infty])\neq 0$.

$\bf{Proof:}$ Assume that $S_{q}$ lies in a tube with rank $d$.

(i) By Riemman-Roch Theorem, $\mbox{Hom}(E, S_{q}[-d])\neq 0$, there has a non-zero morphism
$f': E\to S_{q}[-d]$. Since $\mbox{Ext}^{1}(E, \mathcal{C}^{(q)})=0$, $f'$ can be extended
to a non-zero morphism $f: E\to S_{q}[-\infty]$. Thus $\mbox{Hom}(E, S_{q}[-\infty])\neq 0$.
By Lemma 3.2(ii), $\mbox{Ext}^{1}(E, S_{q}[-\infty])=0$.

(iii) If $\mu{(E)}>q$, $\mbox{Hom}(E, S_{q}[-\infty])=\underleftarrow{\mbox{lim}}\mbox{Hom}(E, S_{q}[-i])=0$.
By Riemman-Roch Theorem, $\mbox{Hom}(S_{q}[-d], \tau E)\neq 0$. So $\mbox{Hom}(S_{q}[-\infty], \tau E)\neq 0$, thus $\mbox{Ext}^{1}(E, S_{q}[-\infty])\neq 0$.

(ii) If $E$ lies in a different tube with $S_{q}$, obviously $\mbox{Hom}(E, S_{q}[-\infty])=\mbox{Ext}^{1}(E, S_{q}[-\infty])$
$=0$. Otherwise $\mbox{Hom}(E, S_{q}[-\infty])=\underleftarrow{\mbox{lim}}\mbox{Hom}(E, S_{q}[-i])=0$.
If $E\neq \tau^{-1}S_{q}$, we get $\mbox{Ext}^{1}(E, S_{q}[-\infty])$
$=0$ since $\mbox{Ext}^{1}(E, S_{q}[-i])=0$ for each $i$. If $E=\tau^{-1}S_{q}$, applying $\mbox{Hom}(-, S_{q})$ to the canonical exact sequence
\begin{equation}
0\rightarrow (\tau S_{q})[-\infty]\rightarrow S_{q}[-\infty]\rightarrow S_{q}\rightarrow 0
\end{equation}
induced by Auslander-Reiten sequences. If $\tau S_{q}\neq S_{q}$, $\mbox{Hom}(S_{q}[-\infty], S_{q})=\mbox{Hom}(S_{q}, S_{q})\neq 0$ and $\mbox{dim}_{k}\mbox{Hom}(S_{q}[-\infty], S_{q})=1$. If not, applying $\mbox{Hom}(S_{q}, -)$ to (3.3), we have the exact sequence
\[\xymatrix{\mbox{Ext}^{1}(S_{q}, S_{q}[-\infty])\ar[r]&\mbox{Ext}^{1}(S_{q}, S_{q}[-\infty])\ar[r]^{g}&\mbox{Ext}^{1}(S_{q}, S_{q})\ar[r]&0\ .}\]
By the similar consideration as Lemma 3.2(ii), $g$ is a monomorphsim. So there has $\mbox{dim}_{k}\mbox{Ext}^{1}(S_{q},$
$S_{q}[-\infty])=1$.$\hfill\square$

Combining the results of Proposition 3.4 and Proposition 3.5, we have

{\bf{Corollary 3.6}} Let $q\in \mathbb{Q}\bigcup\{\infty\}$, then
$$({\bigcap\limits_{S_{q}}}^{\bot}S_{q}[\infty])\cap\mbox{coh}(\mathbb{X})=\mathcal{C}^{(q)}=({\bigcap\limits_{S_{q}}}{}^{\bot}S_{q}[-\infty])\cap
\mbox{coh}(\mathbb{X}),$$
where $S_{q}$ runs through all quasi-simple sheaves of slope $q$, $^{\bot}\mathcal{X}=\{F\in\mbox{Qcoh}\mathbb{X}\ |\ \mbox{Hom}(F, \mathcal{X})$
$=\mbox{Ext}^{1}(F, \mathcal{X})=0\}$.

Moreover, we obtain that

{\bf{Corollary 3.7}} Pr\"{u}fer sheaves and adic sheaves are
indecomposable.

$\bf{Proof:}$ Let $S$ be a quasi-simple sheaf, assume $S[\infty]$ is decomposable, writes $S[\infty]
=U\bigoplus V$ with $U, V\in \mbox{Qcoh}(\mathbb{X})$. By Proposition 3.4 we may assume $U$ satisfies
$\mbox{Hom}(E, U)=0$ for $E\in\mbox{coh}(\mathbb{X})$ with $\mu(E)=q$.
But there exists a surjective morphism from $\bigoplus
S[i]$ to $U$, so it is impossible. Therefore
$S[\infty]$ is indecomposable. Dually adic sheaves are also indecomposable
.$\hfill\square$

{\bf{Corollary 3.8}}  Let $q$, $r\in \mathbb{Q}\bigcup\{\infty\}$,
$S_{q}$, $S'_{r}$ be quasi-simple sheaves of slope $q$ and $r$
respectively. Then

(i) If $q<r$, then $\mbox{Hom}(S_{q}[\infty], S'_{r}[\infty])\neq
0$;

(ii) If $q=r$, then $\mbox{Hom}(S_{q}[\infty], S'_{r}[\infty])\neq
0$ when $S_{q}$, $S'_{r}$ lie in the same tube, otherwise
$\mbox{Hom}(S_{q}[\infty], S'_{r}[\infty])=0$;

(iii) If $q>r$, then $\mbox{Hom}(S_{q}[\infty], S'_{r}[\infty])=0$.

$\bf{Proof:}$ (ii), (iii) is obvious.

(i) By Proposition 3.4, there is a non-zero morphism
from $S_{q}$ to $ S'_{r}[\infty]$. Noticing that for
$i\in\mathbb{N}$, the second rows of the following diagrams always are split,
so there inductively have non-zero maps from
$S_{q}[i+1]$ to $ S_{r}'[\infty]$,
\[\xymatrix{0\ar[r]&S_{q}[i]\ar[r]\ar[d]&S_{q}[i+1]\ar[r]\ar[d]\ar@{.>}[dl]&E_{i}\ar@{=}[d]\ar[r]&0\\
0\ar[r]&S'_{r}[\infty]\ar[r]&H\ar[r]&E_{i}\ar[r]&0.}\]
It implies $\mbox{Hom}(S_{q}[\infty], S_{r}'[\infty])\neq 0$. $\hfill\square$

{\bf{Remark 3.9}} There has the dual property of the morphisms between adic sheaves as Corollary 3.8 which was not showed here.

\section{The Purity of generic, Pr\"{u}fer and
adic sheaves}

Recall that the pure-injective object in a locally
finitely presented category was defined as follows.

{\bf{Definition 4.1}}  (see [6])  Let $\mathcal{A}$ be a locally
finitely presented category, $\mbox{fp}(\mathcal{A})$ be the
subcategory of $\mathcal{A}$ consisting of all finitely presented
objects.

(i) A sequence $0\rightarrow A\rightarrow B\rightarrow C\rightarrow
0$ in $\mathcal{A}$ is called pure-exact if $0\rightarrow
\mbox{Hom}_{\mathcal{A}}(X, A)\rightarrow\mbox{Hom}_{\mathcal{A}}(X,
B)\rightarrow \mbox{Hom}_{\mathcal{A}}(X, C)\rightarrow 0$ is exact
for all $X\in \mbox{fp}(\mathcal{A})$.

(ii) An object $A\in\mathcal{A}$ is called pure-injective if every
pure-exact sequence $0\rightarrow A\rightarrow B\rightarrow
C\rightarrow 0$ is split.

(iii) An object $A\in\mathcal{A}$ is called $\Sigma$-pure-injective
if $\bigoplus_{I} A$ is pure-injective for any set $I$.

Obviously, $\Sigma$-pure-injective object is pure-injective. Noticing that
for a locally finitely presented category $\mathcal{A}$ with
products, the subgroup of finite definition of
$\mbox{Hom}_{\mathcal{A}}(X, A)$ for any $A\in\mathcal{A}$ and $X\in
\mbox{fp}(\mathcal{A})$ is defined as the image of the morphism
$\alpha^{\ast}:\mbox{Hom}_{\mathcal{A}}(Y, A)\rightarrow
\mbox{Hom}_{\mathcal{A}}(X, A)$, arising from a morphism $\alpha:
X\rightarrow Y$ from $X$ to any object $Y\in
\mbox{fp}(\mathcal{A})$. Notice that a subgroup of finite definition
of $\mbox{Hom}(X, A)$ is a sub $\mbox{End}(A)$-module of $\mbox{Hom}(X, A)$.
There has the following property.

 {\bf{Lemma 4.2}}  (see [6]) Let $A$ be an object in a locally finitely presented category
$\mathcal{A}$ with products, then $A$ is $\Sigma$-pure-injective if
and only if $A$ satisfies the descending condition for the subgroup
of finite definition of $\mbox{Hom}_{\mathcal{A}}(X, A)$ for any
$X\in\mbox{fp}(\mathcal{A})$.

Since $\mbox{Qcoh}(\mathbb{X})$ is a locally finitely presented
category with products, we can also consider purity of quasi-coherent
sheaves, and we get

{\bf{Proposition 4.3}} A pure-exact sequence in $\mbox{Qcoh}(\mathbb{X})$ is an
exact sequence.

{\bf{Proof:}} Let $0\rightarrow A\rightarrow B\rightarrow C\rightarrow 0$ be a pure-exact
sequence. By definition, there has an  exact sequence $0\rightarrow\mbox{Hom}(
\oplus L, A)\rightarrow \mbox{Hom}(
\oplus L, B)\rightarrow\mbox{Hom}(
\oplus L, C)\rightarrow 0$ where $L$ runs through all line bundles in
$\mbox{coh}(\mathbb{X})$. It implies $0\rightarrow\Gamma(A)\rightarrow
\Gamma(B)\rightarrow\Gamma(C)\rightarrow 0$ is an exact sequence where $\Gamma(-)$
is the graded global section functor. By sheafication, it finishes the proof.$\hfill\square$

In this section, we discuss the purity of generic sheaves, Pr\"{u}fer sheaves and
adic sheaves. The notion of generic sheaves on $\mathbb{X}$ was introduced in [14], [7]. If $\mathbb{X}$ is a weighted projective line of genus one,
let $T$ be a tilting sheaf on $\mathbb{X}$, by definition, an
indecomposable quasi-coherent sheaf $G$ is called generic if $G$ is
not a coherent sheaf, and $\mbox{Hom}(T, G)$ and $\mbox{Ext}^{1}(T,
G)$ have finite $\mbox{End}(G)$-length; When $\mathbb{X}$ is an elliptic curve, $G$ is called generic if $G$ is not a coherent sheaf, but $\mbox{Hom}(E, G)$ have finite $\mbox{End}(G)$-length for each coherent sheaf $E$.

There exists a unique indecomposable
generic sheaf $G_{q}$ of slope $q$ for each
$q\in\mathbb{Q}\bigcup\{\infty\}$ under isomorphism. In particular,
the sheaf $K$ of rational functions on $\mathbb{X}$ is the generic
sheaf of slope $\infty$, and there exists an automorphism
$\Phi_{q\infty}$ of the bounded derived category
$\mbox{D}^{b}(\mbox{Qcoh}(\mathbb{X}))$ such that
$G_{q}=\Phi_{q\infty}(K)$.
Moreover, the morphisms between generic sheaf $G_{q}$ and coherent sheaves as follows.

{\bf{Lemma 4.4}} (see [14], [7])  Let $E\in \mbox{coh}(\mathbb{X}$) be indecomposable.

(i) If $\mu(E)<q$, then $\mbox{Ext}^{1}(E, G_{q})=0$ and
$\mbox{Hom}(E, G_{q})\neq 0$.

(ii) If $\mu(E)=q$, then $\mbox{Hom}(E, G_{q})=0=\mbox{Ext}^{1}(E,
G_{q})$.

(iii) If $\mu(E)>q$, then $\mbox{Hom}(E, G_{q})=0$ and
$\mbox{Ext}^{1}(E, G_{q})\neq 0$.

{\bf{Theorem 4.5}} Generic sheaves, Pr\"{u}fer sheaves
are $\Sigma$-pure-injective, adic sheaves are pure-injecitve.

$\bf{Proof:} $  Since each generic sheaf $G_{q}$ is of finite length
over End$(G_{q})$ by definition, generic sheaves are $\Sigma$-pure-injective
by Lemma 4.2.

Next, we prove Pr\"{u}fer sheaves are $\Sigma$-pure-injective
sheaves. Let $E$ be an indecomposable coherent sheaf and
$S_{q}[\infty]$ be a Pr\"{u}fer sheaf with quasi-simple sheaf
$S_{q}$ of slope $q$, it suffices to show $\mbox{Hom}(E,
S_{q}[\infty])$ is artinian over $\mbox{End}(S_{q}[\infty])$. If $\mu(E)>q$, $\mbox{Hom}(E, S_{q}[\infty])=0$;
If $\mu(E)=q$, assume $E$ is a quasi-simple sheaf. When $E\neq S_{q}$, $\mbox{Hom}(E, S_{q})=0$. When $E=S_{q}$, $\mbox{dim}_{k}\mbox{Hom}(E, S_{q}[\infty])=1$. For $\mbox{coh}(\mathbb{X})$ is a $k$-linear category, $\mbox{Hom}(E, S_{q}[\infty])$ is artinian over $\mbox{End}(S_{q}[\infty])$. $\mathcal{C}^{(q)}$ is a uniserial subcategory and $\mbox{Ext}^{1}(\mathcal{C}^{(q)}, S_{q}[\infty])=0$, so for each $E$, $\mbox{Hom}(E, S_{q}[\infty])$ is artinian over $\mbox{End}(S_{q}[\infty])$. Next We only need to consider the situation when $\mu(E)<q$ and $S_{q}$ lies
in the tube of rank one.

In this case, there exist an exact commutative diagram
\[\xymatrix{0\ar[r]&S_{q}\ar[r]\ar@{=}[d]&S_{q}[i]\ar^{\phi_{i}}[r]\ar[d]&S_{q}[i-1]\ar[r]\ar[d]&0\\
 0\ar[r]&S_{q}\ar[r]&S_{q}[i+1]\ar^{\phi_{i+1}}[r]&S_{q}[i]\ar[r]&0}\]
where $i\in\mathbb{N}$ and $i>1$. Taking direct limit, we obtain a
new exact sequence
\[\xymatrix{0\ar[r]&S_{q}\ar[r]&S_{q}[\infty]\ar^{\phi}[r]&S_{q}[\infty]\ar[r]&0.}\]
satisfies $\phi|_{S_{q}[i]}=\phi_{i}$.
Since $\mbox{coh}(\mathbb{X})$ is a Krull-Schmidt category, we can
choose a basis of $\mbox{Hom}(E, S_{q})$, write $e_{11}, e_{12}, \ldots,
e_{1n}$.
Applying $\mbox{Hom}(E, -)$ to the exact sequence
\[\xymatrix{0\ar[r]& S_{q}\ar[r]&S_{q}[2]\ar^{\phi_{2}}[r]&S_{q}\ar[r]&0}\]
we obtain a set of elements $e_{21}, e_{22}, \ldots, e_{2n}$ in
$\mbox{Hom}(E, S_{q}[2])$ satisfies $\phi e_{2j}=e_{1j}$.
For $i>2$, inductively, choose elements $e_{i1}, e_{i2}, \ldots, e_{in}$
in $\mbox{Hom}(E, S_{q}[i])$ such that $\phi
e_{ij}=e_{(i-1)j}$. Let $\bigoplus_{i\in\mathbb{N}} kh_{i}$ be a
free $k$-module with basis $(h_{i})_{i\in\mathbb{N}}$ where $h_{i}$ satisfies
$\phi h_{i}=h_{i-1}$ for $i>1$ and $\phi h_{1}=0$. Then it induces a
$k[\phi]$-module structure on $\bigoplus_{i\in\mathbb{N}} kh_{i}$.
Now set $H_{j}=\bigoplus_{i\in\mathbb{N}} kh_{ij}$ be a copy of the
$k[\phi]-$module $\bigoplus_ {i\in\mathbb{N}} kh_{i}$ for every
$j\in\{1, \ldots, n\}$. Then the assignment $h_{ij}\longmapsto
e_{ij}$ induces an epimorphism $\bigoplus^{n}_{j=1} H_{j}\rightarrow
\mbox{Hom}(E, S_{q}[\infty])$ of $k[\phi]$-modules. Since $\bigoplus^{n}_{j=1} H_{j}$ is
artinian over $k[\phi]$ [12], we have $\mbox{Hom}(E, S_{q}[\infty])$
is artinian over $k[\phi]$. So $\mbox{Hom}(E, S_{q}[\infty])$ is
artinian over End$(S_{q}[\infty])$.

At last, we show that each $S_{q}[-\infty]$ is a pure-injective
object. Let
\[\xymatrix{0\ar[r]& S_{q}[-\infty]\ar^{f}[r]&E\ar[r]&
F\ar[r]& 0}\]
be a pure-exact sequence. The pushout of $f$ and the canonical morphism
$S_{q}[-\infty]\rightarrow S_{q}[-i]$ induces a commutative diagram
\[\xymatrix{0\ar[r]&S_{q}[-\infty]\ar[d]\ar^{f}[r]&E\ar[r]\ar[d]&F\ar[r]\ar@{=}[d]&0\\
0\ar[r]&S_{q}[-j]\ar^{f_{j}}[r]\ar[d]&E_{j}\ar[r]\ar[d]&F\ar[r]\ar@{=}[d]&0\\
0\ar[r]&S_{q}[-i]\ar^{f_{i}}[r]&E_{i}\ar[r]&F\ar[r]&0\ .}\]
Obviously, the rows of commutative diagrams are all pure-exact sequences. Since for each $i\in\mathbb{N}, S_{q}[-i]$ is a
pure-injective object, the rows are split, i.e. there exist $h_{i}:
E_{i}\rightarrow S_{q}[-i]$ for $i\in\mathbb{N}$ satisfying $h_{i}f_{i}=\mbox{id}_{S_{q}[-i]}$.
By the universal property of inverse limit, we obtain a morphism $h:E\rightarrow
S_{q}[-\infty]$ satisfying $hf=\mbox{id}_{S_{q}[-\infty]}$. So the first row of the commutative
diagram is split, i.e. $S_{q}[-\infty]$ is a pure-injective object.
$\hfill\square$

\section{Relationship between Pr\"{u}fer, adic and generic sheaves}

In this section, we describe the relationship between
Pr\"{u}fer sheaves, adic sheaves and generic sheaves.

Recall that, for a quasi-coherent sheaf $X$, its torsion part
$tX$ is defined as the sum of all subobjects of $X$ having finite
length. If $tX=0$, i.e. $\mbox{Hom}(S, X)=0$ for each simple sheaf
$S$, then $X$ is called torsion-free. $X$ is called divisible if
$\mbox{Ext}^{1}(S, X)=0$ for each simple sheaf $S$. We extend these
definitions to the following.

{\bf{Definition 5.1}} Let $G\in \mbox{Qcoh}(\mathbb{X})$, $G$ is
called $q$-torsion-free if $\mbox{Hom}(E, G)=0$ for $E\in$
$\mbox{coh}\mathbb{(X)}$ and $\mu(E)\geq q$. $G$ is called
$q$-divisible if $\mbox{Ext}^{1}(E, G)=0$ for $E\in$
$\mbox{coh}\mathbb{(X)}$ and $\mu(E)=q$, i.e. $\mbox{Ext}^{1}(S_{q},
G)=0$ for each quasi-simple sheaf $S_{q}$ of slope $q$.

We having following theorem.

{\bf{Theorem 5.2}} Let $G$ be a $q$-torsion-free divisible sheaf and
$G_{q}$ be the generic sheaf of slope $q$. Then $G=\oplus G_{q}$.

$\bf{Proof:}$ Since $G$ is $q$-torsion-free, we get
$G=\underrightarrow{\mbox{lim}} E_{i}$ for $E_{i}\in
\mbox{coh}\mathbb{(X)}$ with $\mu(E_{i})<q$. Noticing that there is
an automorphism $\Phi_{q\infty}$ of
$\mbox{D}^{b}(\mbox{Qcoh}(\mathbb{X}))$ which sends
$\mathcal{C}^{(\infty)}$ to $\mathcal{C}^{(q)}$ and
$G_{q}=\Phi_{q\infty}(K)$ for $q\in \mathbb{Q}$, we have
$\Phi_{q\infty}^{-1}(G)=\underrightarrow{\mbox{lim}}\Phi_{q\infty}^{-1}(E_{i})\in
\mbox{Qcoh}{\mathbb{(X)}}$.

On the other hand, let $S$ be a simple sheaf, then
$\Phi_{q\infty}(S)$ is a quasi-simple sheaf of slope $q$. And we
have $\mbox{Ext}^{1}(S,
\Phi_{q\infty}^{-1}(G))=\mbox{Ext}^{1}(\Phi_{q\infty}(S), G)=0$ and
$\mbox{Hom}(S, \Phi_{q\infty}^{-1}(G))=\mbox{Hom}(\Phi_{q\infty}(S),
G)=0$. So $\Phi_{q\infty}^{-1}(G)$ is a torsion-free divisible sheaf
of slope $\infty$. Using the similar method as [7], we have a torsion-free divisible sheaf is a direct sum of rational function sheaf, i.e. $\Phi_{q\infty}^{-1}(G)=\oplus K$,
so $G=\oplus\Phi_{q\infty}(K)=\oplus G_{q}$.
$\hfill\square$

 Let $S_{q}$ be a quasi-simple sheaf of slope $q$ which lies in a tube of
 rank $d$ and
$S_{q}[\infty]$ be the corresponding Pr\"{u}fer sheaf. Now we will describe the relationship between Pr\"{u}fer sheaves and generic
sheaves as follows.

{\bf{Theorem 5.3}} There are two exact sequences as follows in
$\mbox{Qcoh}(\mathbb{X})$:
\begin{center}
$0\rightarrow S_{q}\rightarrow S_{q}[\infty]\rightarrow
(\tau^{-1}S_{q})[\infty]\rightarrow 0$ \ \ and  \ \ $0\rightarrow
S_{q}[d]\rightarrow S_{q}[\infty]\rightarrow
S_{q}[\infty]\rightarrow 0$
\end{center}
which produce two inverse systems $\{(\tau^{i}S_{q})[\infty]\ |\ i\in
\mathbb{N}\}$ and $\{S_{q}[\infty]\ |\ i\in \mathbb{N}\}$. Moreover,
we have $\underleftarrow{\mbox{lim}}(\tau^{i}S_{q})[\infty]=\oplus
G_{q}$ and $\underleftarrow{\mbox{lim}}S_{q}[\infty]=\oplus
G_{q}$.

$\bf{Proof:} $ According the Auslander-Reiten quiver of $\mathcal{C}^{(q)}$ we
have an exact commutative diagram
\[\xymatrix{0\ar[r]&S_{q}\ar[r]\ar@{=}[d]&S_{q}[i]\ar[r]\ar[d]&(\tau^{-1}S_{q})[i-1]\ar[r]\ar[d]&0\\
0\ar[r]&S_{q}\ar[r]&S_{q}[i+1]\ar[r]&(\tau^{-1}S_{q})[i]\ar[r]&0}\]
for $i\in\mathbb{N}$. Taking the direct limit, we obtain the
required first exact sequence. We also have another exact
commutative diagram
\[\xymatrix{0\ar[r]&S_{q}[d]\ar[r]\ar@{=}[d]&S_{q}[id]\ar[r]\ar[d]&S_{q}[(i-1)d]\ar[r]\ar[d]&0\\
0\ar[r]&S_{q}[d]\ar[r]&S_{q}[(i+1)d]\ar[r]&S_{q}[id]\ar[r]&0}\] for
$i\in\mathbb{N}$. Taking the direct limit, we obtain the required
second exact sequence.

For the first sequence, we get an inverse system
$\{(\tau^{i}S_{q})[\infty]\ |\ i\in \mathbb{N}\}$, and another inverse
system $\{S_{q}[\infty]\ |\ i\in\mathbb{N}\}$ from the second sequence. By Theorem 5.2,
we only need to show that
$\underleftarrow{\mbox{lim}}(\tau^{i}S_{q})[\infty]$ and
$\underleftarrow{\mbox{lim}}S_{q}[\infty]$ are $q$-torsion-free
divisible sheaves.

By Lemma 3.2 we know that
$\mbox{Ext}^{1}(E, \underleftarrow{\mbox{lim}}S_{q}[\infty])=0$
and $\mbox{Ext}^{1}(E, \underleftarrow{\mbox{lim}}(\tau^{i}S_{q})[\infty])$

\noindent $=0$
for $E$ is a coherent sheaf of slope $q$, i.e. they are $q$-divisible.

Let $E\in \mbox{coh}(\mathbb {X})$ and $\mu(E)\geq q$, then
$\mbox{Hom}(E, S_{q}[\infty])\neq 0$ when $E=S_{q}$. But $S_{q}$ is
a subobject of $S_{q}[d]$, so $\mbox{Hom}(S_{q},
\underleftarrow{\mbox{lim}}S[\infty])=0$, and then
$\underleftarrow{\mbox{lim}}S_{q}[\infty]$ is $q$-torsion-free.

By Proposition 3.4, we have $\mbox{Hom}(E,
\underleftarrow{\mbox{lim}}(\tau^{i}S_{q})[\infty])=0$ when
$\mu(E)\geq q$. So
$\underleftarrow{\mbox{lim}}(\tau^{i}S_{q})[\infty]$ is also
$q$-torsion-free. $\hfill\square$

By duality, we obtain the relationship between adic sheaves and generic
sheaves.

{\bf{Theorem 5.4}} There are two exact sequences as follows in
$\mbox{Qcoh}(\mathbb{X})$:
\begin{center}
$0\rightarrow (\tau S_{q})[-\infty]\rightarrow S_{q}[-\infty]\rightarrow
S_{q}\rightarrow 0$ \ \ and  \ \ $0\rightarrow
S_{q}[-\infty]\rightarrow S_{q}[-\infty]\rightarrow
S_{q}[-d]\rightarrow 0$
\end{center}
which produce two direct systems $\{(\tau^{-i}S_{q})[-\infty]\ |\ i\in
\mathbb{N}\}$ and $\{S_{q}[-\infty]\ |\ i\in \mathbb{N}\}$. Moreover,
we have $\underrightarrow{\mbox{lim}}(\tau^{-i}S_{q})[-\infty]=\oplus
G_{q}$ and $\underrightarrow{\mbox{lim}}S_{q}[-\infty]=\oplus
G_{q}$.

Moreover, there has

{\bf{Corollary 5.5}}   There is an exact sequence
\begin{center}
$0\rightarrow(\tau S_{q})[-\infty]\rightarrow\oplus
G_{q}\rightarrow S_{q}[\infty]\rightarrow 0$.
\end{center}

$\bf{Proof:} $ This proof can also be seen in [3]. Firstly, by
Theorem 5.3 we have an exact sequence
\[\xymatrix{0\ar[r]&S_{q}[d]\ar[r]&S_{q}[\infty]\ar^{\phi}[r]&
S_{q}[\infty]\ar[r]&0}\]
with $\mbox{Ker}\phi^{i}=S_{q}[id]$.
Since $S_{q}[id]=(\tau S_{q})[-id], i\in\mathbb{N}$, we have commutative diagrams
\[\xymatrix{0\ar[r]&(\tau S_{q})[-(i+1)d]\ar[r]\ar[d] & S_{q}[\infty]\ar^{\phi^{i+1}}[r]\ar^{\phi}[d]& S_{q}[\infty]\ar@{=}[d]\ar[r] &0\\
0\ar[r]&(\tau S_{q})[-id]\ar[r] & S_{q}[\infty]\ar^{\phi^{i}}[r]& S_{q}[\infty]\ar[r] &0\ .}\]

\noindent Since $\{(\tau S_{q})[-id]\ |\ i\in\mathbb{N}\}$
satisfies Mittag-Leffler condition, taking inverse limit, by Corollary 4.3 in [11],
we have the required exact sequence.
$\hfill\square$

Next, we consider the morphisms between Pr\"{u}fer sheaves, adic sheaves and
generic sheaves. We can obtain the following results.

{\bf{Corollary 5.6}}  Let $q\in \mathbb{Q}\bigcup\{\infty\}$,
$S_{q}$ be a quasi-simple sheaf of slope $q$,

(i) If $q<r$, then $\mbox{Hom}(S_{q}[\infty], G_{r})\neq 0$ and
$\mbox{Hom}(G_{r}, S_{q}[\infty])=0$.

(ii) If $q\geq r$, then $\mbox{Hom}(S_{q}[\infty], G_{r})=0$ and
$\mbox{Hom}(G_{r}, S_{q}[\infty])\neq 0$.

$\bf{Proof:}$ (i)  If $q<r$, there exist exact commutative diagrams
\[\xymatrix{0\ar[r]&S_{q}[i]\ar[r]\ar[d]&S_{q}[i+1]\ar[r]\ar[d]\ar@{.>}[dl]&E_{i}\ar@{=}[d]\ar[r]&0\\
0\ar[r]&G_{r}\ar[r]&H\ar[r]&E_{i}\ar[r]&0}\] for $i\in\mathbb{N}$.
For by Lemma 4.4, there is a non-zero morphism from $S_{q}$ to $G_{r}$ and
the second rows of commutative diagrams are split, there are
non-zero maps from $S_{q}[i+1]$ to $G_{r}$. It implies
$\mbox{Hom}(S_{q}[\infty], G_{r})\neq 0$. By [14], $G_{q}$ can be
written as a direct limit $G_{q}=\underrightarrow{\mbox{lim}} E_{i}$
where $E_{i}\in\mbox{coh}(\mathbb{X})$ with $\mbox{lim}_{i}
\mu(E_{i})=q$. So $\mbox{Hom}(G_{r},
S_{q}[\infty])=\mbox{Hom}(\underrightarrow{\mbox{lim}} E_{i},
S_{q}[\infty])=\underleftarrow{\mbox{lim}}\mbox{Hom}(E_{i},
S_{q}[\infty])=0$.

(ii) If $q=r$, using Theorem 5.3, the result is obvious.
If $q>r$, by Lemma 4.4, we know that $\mbox{Hom}(S_{q}[\infty],
G_{r})=0$. And $\mbox{Hom}(G_{r},
S_{q})=\mbox{DExt}^{1}(\tau^{-1}S_{q}, G_{r})\neq 0$, so
$\mbox{Hom}(G_{r}, S_{q}[\infty])\neq 0$. $\hfill\square$

{\bf{Corollary 5.7}} Let $q\in \mathbb{Q}\bigcup\{\infty\}$,
$S_{q}$ be a quasi-simple sheaf of slope $q$,

(i) If $q\leq r$, then $\mbox{Hom}(S_{q}[-\infty], G_{r})\neq 0$ and
$\mbox{Hom}(G_{r}, S_{q}[-\infty])=0$.

(ii) If $q>r$, then $\mbox{Hom}(S_{q}[-\infty], G_{r})=0$ and
$\mbox{Hom}(G_{r}, S_{q}[-\infty])\neq 0$.

{\bf{Proof:}} (i) If $q<r$, there exists a non-zero morphism from $S_{q}$ to $G_{r}$, and
since there has a canonical surjective morphism from $S_{q}[-\infty]$ to $S_{q}$, so
$\mbox{Hom}(S_{q}[-\infty], G_{r})\neq 0$. By Lemma 4.4 and Lemma 3.2,  $\mbox{Hom}(G_{r}, S_{q}[-\infty])=0$. If $q=r$, by Theorem 5.4 $S_{q}[-\infty]$ is a subobject of direct sum of $G_{q}$
, i.e. $\mbox{Hom}(S_{q}[-\infty], G_{q})\neq 0$. Obviously, $\mbox{Hom}(G_{q}, S_{q}[-\infty])=0$.

(ii) If $q>r$, there exists a non-zero morphism from $G_{r}$ to $S_{q}$. Since $\mbox{Ext}^{1}
(G_{r}, \mathcal{C}^{(q)})=0$, it can be extended to a non-zero morphism from $G_{r}$ to
$S_{q}[-\infty]$. Therefore, we obtain that $\mbox{Hom}(G_{r}, S_{q}[-\infty])\neq 0$.
By Proposition 3.5, $S_{q}[-\infty]$ can be written a direct limit of coherent sheaves of slope
greater than $r$, so $\mbox{Hom}(S_{q}[-\infty], G_{r})=0$.$\hfill\square$

\section{The construction of generic sheaves}

In this section, we always assume that $q\in\mathbb{Q}\bigcup\{\infty\}$. For each $q$, denote

$\mathcal{C}_{q}=
\{F\in\mbox{Qcoh}(\mathbb{X})|F$ is $q^{\prime}$-torsion-free where $q^{\prime}\in \mathbb{Q}\bigcup \{\infty\}$ and $q^{\prime}>q \}$,

$\mathcal{Q}_{q}=\{F\in\mbox{Qcoh}(\mathbb{X})| F$ is a factor of direct sum of sheaves in $\bigcup_{q^{\prime}\in \mathbb{Q}\bigcup \{\infty\},q^{\prime}>q}\mathcal{C}^{(q^{\prime})}\}$,

$w_{q}=
\{W\in\mbox{Qcoh}(\mathbb{X})|W\in\mathcal{C}_{q}$ and $W$ is $q$-divisible$\}$.

\noindent And let $w_{q}'\subseteq w_{q}$ be the full subcategory of all direct sums of Pr\"{u}fer sheaves of slope $q$.

In this section, we will show the $w_{q}$-approximation of each quasi-coherent sheaf, and then provide two methods to construct generic sheaves over a weighted
projective line or an elliptic curve by using coherent sheaves and Pr\"{u}fer sheaves. Firstly, we obtain
the following important property.

{\bf{Proposition 6.1}} $(\mathcal{Q}_{q}, \mathcal{C}_{q})$ is a split torsion pair.

{\bf{Proof:}} Notice that when $q=\infty$, we have $\mathcal{C}_{q}=\mbox{coh}(\mathbb{X})$ and $\mathcal{Q}_{q}=0$.
By [3], $(\mathcal{Q}_{q}, \mathcal{C}_{q})$ is always a torsion pair. We only need consider the cases $q\in\mathbb{Q}$.

Let $\eta:0 \rightarrow F\rightarrow G\rightarrow H\rightarrow 0$ be
an exact sequence with $F\in\mathcal{Q}_{q}$ and
$H\in\mathcal{C}_{q}$. We first show if $H$ is a coherent sheaf,
then $\eta$ is split. Without loss of generality, we assume $H$ is
indecomposable. If $\eta$ is not split, there exists a commutative
diagram
\[\xymatrix{0\ar[r]&F\ar[r]\ar^{\alpha}@{.>}[d]&G\ar[r]\ar@{.>}[d]&H\ar[r]\ar@{=}[d]&0\\
0\ar[r]&\tau H\ar[r]&E\ar[r]&H\ar[r]&0}\]
where the second row is an Auslander-Reiten sequence. Since $\tau H\in\mathcal{C}_{q}$, then $\alpha=0$,
which is impossible. So $\eta$ is split. Secondly, we assume that $H$ is a quasi-coherent sheaf, we can
write $H=\underrightarrow{\mbox{lim}}H_{i}$, where $H_{i}$ is the coherent subsheaf of $H$. Certainly, the slope of each $H_{i}$ is not greater than
$q$. So by Lemma 3.2, we obtain $\mbox{Ext}^{1}(H, F)=\mbox{Ext}^{1}(\underrightarrow{\mbox{lim}}H_{i}, F)=0$.
Therefore $(Q_{q}, \mathcal{C}_{q})$ is a split torsion pair.$\hfill\square$

{\bf{Remark 6.2}} Similarly, let $\mathcal{C}_{q}'=
\{F\in\mbox{Qcoh}(\mathbb{X})|F$ is $q$-torsion-free\} and $\mathcal{Q}_{q}'=\{F\in\mbox{Qcoh}(\mathbb{X})| F$ is a factor of direct sum of sheaves in $\bigcup_{q^{\prime}\in \mathbb{Q}\bigcup \{\infty\}, q^{\prime}\geq q}\mathcal{C}^{(q^{\prime})}\}$,
then $(\mathcal{Q}_{q}', \mathcal{C}_{q}')$ is also a split torsion pair.

For any class $\mathcal{Z}$ of quasi-coherent sheaves, we denote by $l(\mathcal{Z})$ the class of all
quasi-coherent sheaves $F$ with $\mbox{Hom}(F, \mathcal{Z})=0$ and $r(\mathcal{Z})$ the class of all
quasi-coherent sheaves $F$ with $\mbox{Hom}(\mathcal{Z}, F)=0$. Then we have
$(lr(\mathcal{Z}), r(\mathcal{Z}))$
is a torsion pair in $\mbox{Qcoh}(\mathbb{X})$. Moreover, denote by $g(\mathcal{C}^{(q)})$ the
class of all quasi-coherent sheaves $F$ generated by $\mathcal{C}^{(q)}$, that is, $F$ is a factor of direct sums of objects in $\mathcal{C}^{(q)}$. To study the structure of quasi-coherent sheaves in $w_{q}$, we need the following lemma.

{\bf Lemma 6.3} $g(\mathcal{C}^{(q)})=lr(\mathcal{C}^{(q)})$ and then $(g(\mathcal{C}^{(q)}), r(\mathcal{C}^{(q)}))$ is a torsion pair in $\mbox{Qcoh}(\mathbb{X})$.

{\bf Proof:}  To show $g(\mathcal{C}^{(q)})=lr(\mathcal{C}^{(q)})$, we only need to prove that $g(\mathcal{C}^{(q)})$
is closed under extension.

Let $0\rightarrow M\rightarrow N\rightarrow N/M\rightarrow 0$ be
an exact sequence in $\mbox{Qcoh}(\mathbb{X})$ with $M$ and $N/M$ in $g(\mathcal{C}^{(q)})$.
If $N\in\mbox{coh}(\mathbb{X})$, then there exist surjective morphisms $f : E\rightarrow M$ and $g : F\rightarrow N/M$
where $E$ and $F$ lie in $\mathcal{C}^{(q)}$. Since $\mbox{Qcoh}(\mathbb{X})$ is a hereditary category, we have the following commutative diagram
\[\xymatrix{0\ar[r]&E\ar[r]\ar^{f}[d]&H\ar[r]\ar^{f'}[d]&F\ar[r]\ar@{=}[d]&0\\
0\ar[r]&M\ar[r]\ar@{=}[d]&H'\ar[r]\ar^{g'}[d]&F\ar[r]\ar^{g}[d]&0\\
0\ar[r]&M\ar[r]&N\ar[r]&N/M\ar[r]&0\ .}\]
Thus $H\in \mathcal{C}^{(q)}$ and  $g'f'$ is surjective which implies $N\in g(\mathcal{C}^{(q)})$.
Now assume $N$ be a quasi-coherent sheaf. Since $N/M\in g(\mathcal{C}^{(q)})$, there exists a surjective
morphism $h: \bigoplus E_{i}\rightarrow N/M$ with $E_{i}\in \mathcal{C}^{(q)}$. Then $N/M=\bigcup h(E_{i}),$ that is,
we can write $N/M=\bigcup(N_{i}/M)$, where $M\subset N_{i}\subset N$
and $N_{i}/M$ is finitely generated by $\mathcal{C}^{(q)}$. Therefore, we have the following commutative diagram
\[\xymatrix{0\ar[r]&M\ar[r]\ar@{=}[d]&N_{i}\ar[r]\ar[d]&N_{i}/M\ar[r]\ar[d]&0\\
0\ar[r]&M\ar[r]\ar@{=}[d]&\bigcup N_{i}\ar[r]\ar[d]&(\bigcup N_{i})/M\ar[r]\ar[d]^{\iota}&0\\
0\ar[r]&M\ar[r]&N\ar[r]&N/M\ar[r]&0\ .}\]
Since $N/M=\bigcup (N_{i}/M)$, we have $\iota$ is an isomorphism which implies $N=\bigcup N_{i}$.
Thus without loss of generality, we only need to show that if $N/M$ is finitely generated by $\mathcal{C}^{(q)}$ and
$M\in g(\mathcal{C}^{(q)}),$ then $N\in g(\mathcal{C}^{(q)})$. Write $N=\bigcup N_{i},$ where $
\{N_{i}\}$ is a set of filtered subcoherent sheaves of $N$. Thus $N/M=(\bigcup N_{i})/M=
\bigcup (N_{i}/M)=\bigcup (N_{i}\bigcup M)/M$, where the second equality is according to Proposition 11.2 in [16]. Notice that $(N_{i}\bigcup M)/M\in\mbox{coh}(\mathbb{X})$
 and $(N_{i}\bigcup M)/M\subset(N_{i'}\bigcup M)/M$ when $i<i'.$  We get $\mbox{deg}(N_{i'}\bigcup M)/M>
 \mbox{deg}(N_{i}\bigcup M)/M$ or $\mbox{rk}(N_{i'}\bigcup M)/M>
 \mbox{rk}(N_{i}\bigcup M)/M$. Thus there exists $i$ such that $N/M=(N_{i}\bigcup M)/M$, which implies
 $N=N_{i}\bigcup M$. Now write $M=\bigcup M_{j},$ where $
\{M_{j}\}$ is a set of filtered subcoherent sheaves of $M$. There exists $j$ such that $
M_{j}\bigcap N_{i}=M\bigcap N_{i}$. So $(N_{i}\bigcup M_{j^{'}})/N_{i}=N_{i}/(M_{j^{'}}\bigcap N_{i})
=N_{i}/(M\bigcap N_{i})=(N_{i}\bigcup M)/M=N/M$ for $j^{'}\geq j$. We obtain that $N_{i}\bigcup M_{j^{'}}\in g({\mathcal{C}^{(q)}})$. Thus $N\in g(\mathcal{C}^{(q)})$.$\hfill\square$

{\bf{Theorem 6.4}} Each $W\in w_{q}$ is direct sums of Pr\"{u}fer sheaves and the generic sheaf $G_{q}$.

{\bf{Proof:}} Let $W\in w_{q}$, according to Lemma 6.3, there has an exact sequence
\begin{equation}
0\rightarrow W_{1}\rightarrow W\rightarrow W_{2}\rightarrow 0,
\end{equation}
where $W_{1}\in g(\mathcal{C}^{(q)})$ and $W_{2}\in r(\mathcal{C}^{(q)})$. Now, $W_{2}$ is $q$-torsion-free divisible implies $W_{2}=\oplus G_{q}$. Since $W_{1}$ is a direct limit of its subsheaves which lie in $\mathcal{C}^{(q)}$, there exist a quasi-simple sheaf $S_{q}$ and a non-zero morphism $\alpha: S_{q}\rightarrow W_{1}$. Obviously, $\alpha$ must be a monomorphism. If not, $\mbox{Ker}\alpha\neq 0$, $\mbox{Ker}\alpha$ ia a line bundle, so $\mbox{Im}(\alpha)\in\mbox{coh}_{0}(\mathbb{X})$, it is a contradiction with $W_{1}\in\mathcal{C}_{q}$. For each quasi-simple sheaf $S_{q}$, applying $\mbox{Hom}(S_{q}, -)$ to (6.1), there exists a long exact sequence
\[\xymatrix{\mbox{Hom}(S_{q}, W)\ar[r]&\mbox{Hom}(S_{q}, W_{2})\ar[r]&\mbox{Ext}^{1}(S_{q}, W_{1})\ar[r]&\mbox{Ext}^{1}(S_{q}, W)\ar[r]&\mbox{Ext}^{1}(S_{q}, W_{2}).}\]
$W_{2}\in r(\mathcal{C}^{(q)})$ and $W$ is $q$-divisible, thus $\mbox{Ext}^{1}(S_{q}, W_{1})=0$, i.e. $W_{1}$ is $q$-divisible. Taking the pushout of $\alpha$ and the canonical monomorphism $S_{q}\rightarrow S_{q}[2]$, we have the following commutative diagram
\[\xymatrix{0\ar[r]&S_{q}\ar^{\eta}[r]\ar^{\alpha}[d]&S_{q}[2]\ar[r]\ar^{\gamma}[d]&\tau^{-1}S_{q}\ar[r]\ar@{=}[d]&0\\
0\ar[r]&W_{1}\ar^{f}[r]&W_{1}'\ar[r]&\tau^{-1}S_{q}\ar[r]&0\ .}\]
For $W_{1}$ is $q$-divisible, there exists $f': W_{1}'\rightarrow W_{1}$ satisfying $f'f=\mbox{id}_{W_{1}}$. Denoted by $\alpha_{2}=f'\gamma$, then $\alpha=\alpha_{2}\eta$. Moreover, $\alpha_{2}$ is also a monomorphism. In fact, if $\mbox{Ker}\alpha_{2}\neq 0$, $\mbox{Ker}\alpha_{2}$ is $S_{q}$ or a vector bundle with slope less than $q$. If $\mbox{Ker}\alpha_{2}=S_{q}$, it is a contradiction with $\alpha\neq 0$; If $\mbox{Ker}\alpha_{2}$ is a vector bundle with slope less than $q$, $\mbox{Im}\alpha_{2}$ is a sheaf with slope larger than $q$, it is a contradiction with $W_{1}\in\mathcal{C}_{q}$. Thus, $\alpha$ can be extended to a new monomorphism $\alpha_{2}: S_{q}[2]\rightarrow W_{1}$. By induction, given each $i$, $\alpha$ can be extended to the monomorphism $\alpha_{i}: S_{q}[i]\rightarrow W_{1}$. By the universality property of direct limit, there is a morphism $\alpha': S_{q}[\infty]\rightarrow W_{1}$. We claim that $\alpha'$ is a monomorphism.

For each $g: F\rightarrow S_{q}[\infty]$ where $\alpha'g=0$, we show that $g=0$. If $F$ is a coherent sheaf, then $F$ is a finitely presented object, there exists $i$ and a morphism $g_{i}: F\rightarrow S_{q}[i]$ satisfying $\beta_{i}g_{i}=g$.
\[\xymatrix{&F\ar@/_/_{g_{i}}[ld]\ar@/^/^{g_{j}}[rd]\ar^{g}[dd]&\\
S_{q}[i]\ar@/_/_{\alpha_{i}}[ddr]\ar@/_/^{\beta_{i}}[dr]\ar[rr]&&S_{q}[j]\ar@/^/^{\alpha_{j}}[ddl]\ar@/^/_{\beta_{j}}[ld]\\
&S_{q}[\infty]\ar^{\alpha'}[d]&\\
&W_{1}&}\]
According to the commutative diagram, $\alpha_{i}g_{i}=\alpha'\beta_{i}g_{i}=\alpha'g=0$. On the other hand $\alpha_{i}$ is a monomorphism, $g_{i}=0$, then $g=0$. If $F$ is a quasi-coherent sheaf, denoted by $F=\underrightarrow{\mbox{lim}}F_{i}$, where $F_{i}$ is the coherent subsheaf of $F$, $\phi_{i}: F_{i}\rightarrow F$ is the canonical inclusion. Since $\alpha'g\phi_{i}=0$, we know that for each $i$, $g\phi_{i}=0$. By the universality property of direct limit, $g=0$.

In summary, $\alpha$ can be extended to new monomorphism $\alpha': S_{q}[\infty]\rightarrow W_{1}$. On the other hand $\mbox{Coker}\alpha'$ is a direct limit of coherent subsheaves in $\mathcal{C}^{(q)}$, by Lemma 3.2, $\alpha'$ is split, i.e. $S_{q}[\infty]$ is a direct summand of $W_{1}$. Finally we prove that $W_{1}$ is direct sums of Pr\"{u}fer sheaves with slope $q$.

By transfinite induction, we need to construct the subsheaf $V_{i}$ of $W_{1}$, satisfies: $V_{i}$ is a direct limit of coherent subsheaves in $\mathcal{C}^{(q)}$, and $V_{i}$ is $q$-divisible; for any ordinal $i$, $V_{i}=V_{i-1}\oplus S_{q}[\infty]$, where $S_{q}$ is a quasi-simple sheaf in $\mathcal{C}^{(q)}$; for any limit ordinal $i$, $V_{i}=\bigcup\limits_{j<i} V_{j}$. Then $W_{1}=V_{i}\oplus U_{i}$, where $U_{i}$ is also a direct limit of coherent subsheaves in $\mathcal{C}^{(q)}$ and $q$-divisible, the construction will stop when $U_{i}=0$.

Let $V_{0}=0$. Assume $V_{i}$ has been defined with $W_{1}=V_{i}\oplus U_{i}$, for $U_{i}$ is a direct limit of coherent subsheaves in $\mathcal{C}^{(q)}$ and $q$-divisible, by the similar analysis to $W_{1}$, $U_{i}$ contains some direct summand $S_{q}[\infty]$, denote $U_{i}=S_{q}[\infty]\oplus U_{i+1}$ and $V_{i+1}=V_{i}\oplus S_{q}[\infty]$, then $W_{1}=V_{i+1}\oplus U_{i+1}$. for any limit ordinal $i$, assume for any $j$, $j<i$, $V_{j}$ has been defined, let $V_{i}=\bigcup\limits_{j<i}V_{j}$. Obviously $V_{i}$ is also a direct limit of coherent sheaves in $\mathcal{C}^{(q)}$. For any quasi-simple sheaf $S_{q}$, if $\mbox{Ext}^{1}(S_{q}, V_{i})\neq 0$, then $\mbox{Hom}(S_{q}, V_{i}[1])\neq 0$, by the definition of $V_{i}$, there exists $j$, $j<i$, satisfying $\mbox{Hom}(S_{q}, V_{j}[1])\neq 0$, it is a contradiction with $\mbox{Ext}^{1}(S_{q}, V_{j})=0$, i.e. $V_{i}$ is $q$-divisible. By construction $V_{i}$ is direct sums of Pr\"{u}fer sheaves with slope $q$. By transfinite induction, $W_{1}$ is also direct sums of Pr\"{u}fer sheaves with slope $q$.

Since $\mbox{Ext}^{1}(G_{q}, S_{q}[\infty])=0$, (6.1) is split, then $W=W_{1}\oplus W_{2}$, $W$ is direct sums of Pr\"{u}fer sheaves and the generic sheaf $G_{q}$.$\hfill\square$

{\bf{Proposition 6.5}} For $q\in\mathbb{Q}$, each $F\in \mathcal{Q}_{q}$ is generated by $w_{q}'$, i.e. $F$ is a factor of an object in $w_{q}'$.

{\bf{Proof:}} We only need to prove that if $F$ is a coherent
sheaf of slope greater than $q$, $F$ is generated by $w_{q}'$. Denote by $tF$ be the union of all images of non-zero morphism from $\mathcal{C}^{(q)}$ to $F$. If $tF\neq F$, then $F/tF$ is a non-zero coherent sheaf of slope greater than $q$. By Riemman-Roch formula, there is a non-zero morphism from $\mathcal{C}_{q}$ to $
F/tF$, it is a contradiction by Lemma 6.3. So $F$ is generated by $\mathcal{C}^{(q)}$. Since $\mbox{Ext}^{1}(\mathcal{C}^{(q)}, F)=0$, any morphism from $\mathcal{C}^{(q)}$ to $F$ can be extended
to the morphism from Pr\"{u}fer sheaves of slope $q$ to $F$, so $F$ is generated by $w_{q}'$. $\hfill\square$

When considering the left $w_{q}$-approximation of $\mathcal{C}_{q}$, we have the following theorem which can be immediately obtained from Theorem 4.1 in [20] by using similar method.

{\bf{Theorem 6.6}} For each $F\in\mathcal{C}_{q}$, there exists an exact sequence
\[\xymatrix{0\ar[r]&F\ar^{\mu_{q}}[r]&F_{w_{q}}\ar[r]&\oplus_{S_{q}}\oplus S_{q}[\infty]\ar[r]&0}\]
where $\mu_{q}$ is the minimal left $w_{q}$-approximation. Moreover, if $F$ is $q$-torsion-free,
then $F_{w_{q}}$ is also $q$-torsion-free.$\hfill\square$

Denote $\mbox{lg}$ be the length of a module, we obtain the connection between Pr\"{u}fer sheaf and generic sheaf
by the exact sequence as follows.

{\bf{Theorem 6.7}} Let $F\in\mathcal{C}_{q}$ and $F$ is $q$-torsion-free. Moreover, if
\begin{center}
$\mbox{Hom}(F, \oplus G_{q})=\oplus\mbox{Hom}(F, G_{q})$\ and\ $\mbox{lg}_{\mbox{End}(G_{q})}(\mbox{Hom}(F, G_{q}))=n\in\mathbb{Z},$
\end{center}
then there exists an exact sequence
\[\xymatrix{
0\ar[r]&F\ar^{\mu_{q}}[r]& \oplus_{n}G_{q}\ar[r]&\oplus_{S_{q}}\oplus_{e_{S_{q}F}}S_{q}[\infty]\ar[r]&0
}\]
where $e_{S_{q}F}=\mbox{lg (Ext}^{1}(S_{q}, F))_{\mbox{End}(S_{q})}$
and $S_{q}$ runs through all quasi-simple sheaves of slope $q$.

{\bf{Proof:}} By Theorem 6.6, if $F$ is $q$-torsion-free, $F_{w_{q}}$ is direct sums of
$G_{q}$ and $\mu_{q}$ is naturally a minimal $\mathcal{G}$-approximation where $\mathcal{G}
=\{F\in\mbox{Qcoh}(\mathbb{X})\ |\ F$ is direct sums of $G_{q}\}$. On the other hand, if
$\mbox{lg}_{\mbox{End}(G_{q})}(\mbox{Hom}(F, G_{q}))=n$, write $\{e_{1}, e_{2},\ldots, e_{n}\}$
be a basis of $\mbox{Hom}(F, G_{q})$ over $\mbox{End}(G_{q})$. Let $f=(e_{1}, e_{2},\ldots, e_{n})^{\intercal}: F\to \oplus_{n}G_{q}$, for $\mbox{Hom}(F, \oplus G_{q})=\oplus\mbox{Hom}(F, G_{q})$, obviously $f$ is also a minimal $\mathcal{G}$-approximation which implies
$f=\mu_{q}$. This finishes our proof.
$\hfill\square$

According to [14], [7], there exists a linear form over $\mbox{coh}(\mathbb{X})$ as follows:
\begin{center}
$\langle[E],
[G_{q}]\rangle=\mbox{lg}_{\mbox{End}(G_{q})}(\mbox{Hom}(E, G_{q}))
-\mbox{lg}_{\mbox{End}(G_{q})}(\mbox{Ext}^{1}(E, G_{q})),$
\end{center}
where $E\in\mbox{coh}(\mathbb{X})$, $q=d/r$, $d$ and $r$ are coprime
integers, $r\geq 0$, $[G_{q}]=ru+dw$. Notice that, the definition of $u$ and $w$ are different for a weighted projective line of genus one and an elliptic curve.

{\bf Corollary 6.8} Let $F$ be a coherent sheaf
with slope less than $q$. If $F$ satisfies the equation
$d\mbox{rk}(F)-r\mbox{deg}(F)=1$, then there has an exact sequence
\[\xymatrix{
0\ar[r]&F\ar[r]&G_{q}\ar[r]&\oplus_{S_{q}}\oplus_{e_{S_{q}F}}S_{q}[\infty]\ar[r]&0,
}\] where $e_{S_{q}F}=\mbox{lg (Ext}^{1}(S_{q},
F))_{\mbox{End}(S_{q})}$ and $S_{q}$ runs through all quasi-simple
sheaves of slope $q$.

{\bf Proof:} If $F$ is a coherent sheaf with slope less than $q$, $F$ is naturally $q$-torsion-free. Moreover, $\langle [F], [G_{q}]\rangle=\langle [F], ru+dw\rangle=d\mbox{rk}(F)-r\mbox{deg}(F)$. According to Theorem 6.7, we immediately get the required statement.$\hfill\square$

{\bf Remark 6.9} Corollary 6.8 provide a method to construct generic
sheaves. Notice that if $q=\infty$, then $G_{\infty}=K$ and a
coherent sheaf $F$ satisfies the conditions in Corollary 6.8 must be
a line bundle. Corollary 6.8 in fact popularizes the method of
construction for $K$ over elliptic curves in [8].

Next we consider the minimal right $w_{q}-$approximation of $\mathcal{Q}_{q}$.

{\bf{Theorem 6.10}} Let $q\in\mathbb{Q}$. For each $F\in\mathcal{Q}_{q}$, there
exists an exact sequence
\[\xymatrix{0\ar[r]&\oplus G_{q}\ar[r]&\oplus_{S_{q}}\oplus S_{q}[\infty]\ar^{ \ \ \  \ \ \ \beta_{q}}[r]&F\ar[r]&0}\]
where $\beta_{q}$ is the minimal right $w_{q}$-approximation.

{\bf{Proof:}} We only need to prove that there is an exact sequence $0\rightarrow G\rightarrow \oplus_{S_{q}}\oplus S_{q}[\infty]\rightarrow F\rightarrow 0$ satisfying $G$ is $q$-torsion-free, the rest is similar to the proof of Theorem 7.1 in [20].

According to Proposition 6.5, we obtain an exact sequence $0\rightarrow G\rightarrow\oplus_{S_{q}}\oplus S_{q}[\infty]\rightarrow F\rightarrow 0$ where $G$ lies in $\mathcal{C}_{q}$.
Moreover, $G$ can be chosen to be $q$-torsion-free. In fact, if it is not, let $tG$ be the union of all images of non-zero morphism from $\mathcal{C}^{(q)}$ to $G$, then $G/tG$ is $q$-torsion-free. We have a new exact sequence $0\rightarrow G/tG\rightarrow (\oplus_{S_{q}}\oplus S_{q}[\infty])/tG\rightarrow F\rightarrow 0$. Since $tG$ lies in $\mathcal{C}_{q}$ and $tG$ is generated by $\mathcal{C}^{(q)}$, we can write $tG=\underrightarrow{\mbox{lim}} E_{i}$ where $E_{i}$ lies in $\mathcal{C}^{(q)}$, so $(\oplus_{S_{q}}\oplus S_{q}[\infty])/tG$ is $q$-divisible and generated by $\mathcal{C}^{(q)}$ which means  $(\oplus_{S_{q}}\oplus S_{q}[\infty])/tG$ is a direct sum of Pr\"{u}fer sheaves. Therefore, we obtain an exact sequence $0\rightarrow G\rightarrow \oplus_{S_{q}}\oplus S_{q}[\infty]\rightarrow F\rightarrow 0$ where $G$ is $q$-torsion-free and lies in $\mathcal{C}_{q}$. $\hfill\square$

{\bf{Remark 6.11}} Compared to Theorem 6.7, $q$ can not be $\infty$ in Theorem 6.10.

{\bf{Theorem 6.12}} Let $F\in\mathcal{Q}_{q}$, if $F$ satisfies $\mbox{lg}_{\mbox{End}(G_{q})}(\mbox{Ext}^{1}(F, G_{q}))=n$, then there has an exact sequence
\[\xymatrix{0\ar[r]&\oplus_{n} G_{q}\ar[r]&\oplus_{S_{q}}\oplus S_{q}[\infty]\ar^{\ \ \ \ \ \beta_{q}}[r]&F\ar[r]&0,}\]
where $\beta_{q}$ is a minimal right $\omega_{q}$-approximation.

{\bf{Proof:}} By Theorem 6.10, there has an exact sequence
\[\xymatrix{0\ar[r]&\oplus_{I} G_{q}\ar[r]&\oplus_{S_{q}}\oplus S_{q}[\infty]\ar^{\ \ \ \ \beta_{q}}[r]&F\ar[r]&0,}\]
where $\beta_{q}$ is a minimal right $\omega_{q}$-approximation, $I$ is an index set. Applying $\mbox{Hom}(-, G_{q})$, we get a long exact sequence
$$\mbox{Hom}(\oplus_{S_{q}}\oplus S_{q}[\infty], G_{q})\rightarrow\mbox{Hom}(\oplus_{I} G_{q}, G_{q})\stackrel{f}{\rightarrow}\mbox{Ext}^{1}(F, G_{q})\rightarrow\mbox{Ext}^{1}(\oplus_{S_{q}}\oplus S_{q}[\infty], G_{q}),$$
$\mbox{Ext}^{1}(\oplus_{S_{q}}\oplus S_{q}[\infty], G_{q})=0$ and $\mbox{Hom}(\oplus_{S_{q}}\oplus S_{q}[\infty], G_{q})=0$, so $f$ is an isomorphism. On the other hand, by $\mbox{lg}_{\mbox{End}(G_{q})}(\mbox{Ext}^{1}(F, G_{q}))=n$, denote $\{\xi_{1}, \xi_{2}, \ldots, \xi_{n}\}$ as the basis of $\mbox{Ext}^{1}(F, G_{q})$ as $\mbox{End}(G_{q})$-vector space. For $\mbox{Ext}^{1}(F, \prod_{n} G_{q})=\prod_{n}\mbox{Ext}^{1}(F, G_{q})$, $(\xi_{1}, \xi_{2}, \ldots, \xi_{n})$ uniquely determines an exact sequence as follows
\[\xymatrix{\xi:&0\ar[r]&\prod_{n} G_{q}\ar[r]&H\ar[r]&F\ar[r]&0,}\]
Moreover, for each $i$, the canonical morphism $\eta_{i}: \prod_{n} G_{q}\rightarrow G_{q}$ gives the following commutative diagram
\[\xymatrix{\xi:&0\ar[r]&\prod_{n} G_{q}\ar^{\eta_{i}}[d]\ar[r]&H\ar[r]\ar[d]&F\ar[r]\ar@{=}[d]&0\\
\xi_{i}:&0\ar[r]&G_{q}\ar[r]&H_{i}\ar[r]&F\ar[r]&0\ .}\]
For any $\varepsilon$ in $\mbox{Ext}^{1}(F, G_{q})$, there exist $g_{i}: G_{q}\rightarrow G_{q}$ satisfying $\varepsilon=\sum\limits_{i=1}^{n}g_{i}\xi_{i}$, and the commutative diagram
\[\xymatrix{\xi:&0\ar[r]&\prod_{n} G_{q}\ar^{(g_{i})}[d]\ar[r]&H\ar[r]\ar[d]&F\ar[r]\ar@{=}[d]&0\\
\varepsilon:&0\ar[r]&G_{q}\ar[r]&H'\ar[r]&F\ar[r]&0\ .}\]
In particular, If $\varepsilon=0$, then $(g_{i})=0$. Thus, applying $\mbox{Hom}(-, G_{q})$ to $\xi$, $\mbox{Hom}(\prod_{n} G_{q}, G_{q})\rightarrow\mbox{Ext}^{1}(F, G_{q})$ is an isomorphism. In summary, $\mbox{Hom}(\oplus_{I} G_{q},$ 

\noindent $G_{q})=\mbox{Hom}(\oplus_{n} G_{q}, G_{q})$, then $|I|$ equals $n$.$\hfill\square$

{\bf{Theorem 6.13}} Given any coherent sheaf $E$ in $\mathcal{Q}_{q}$, if $r\mbox{deg}(E)-d\mbox{rk}(E)=1$, then there has an exact sequence
\[\xymatrix{0\ar[r]&G_{q}\ar[r]&\oplus_{S_{q}}\oplus S_{q}[\infty]\ar^{\ \ \ \ \ \beta_{q}}[r]&E\ar[r]&0,}\]
where $\beta_{q}$ is a minimal right $\omega_{q}$-approximation.

{\bf{Proof:}} $E$ is a coherent sheaf in $\mathcal{Q}_{q}$, so $\mbox{Hom}(E, G_{q})=0$. we compute $\langle E, G_{q}\rangle$, then
$$\langle E, G_{q}\rangle=-r\mbox{deg}(E)+d\mbox{rk}(E)=-\mbox{lg}_{\mbox{End}(G_{q})}\mbox{Ext}^{1}(E, G_{q}).$$
By Theorem 6.12, when $r\mbox{deg}(E)-d\mbox{rk}(E)=1$, the exact sequence in Theorem 6.13 holds.$\hfill\square$

The theory of infinitely generated tilting module is a focus of area of algebraic representation theory in recent years. In view of the relationship between quasi-coherent sheaves and infinitely generated modules, we give the concept of large tilting sheaves in $\mathbb{X}$.

{\bf{Definition 6.14}} Let $\mathbb{X}$ be a weighted projective line of genus one or an elliptic curve, $T\in\mbox{Qcoh}(\mathbb{X})$. If $\mbox{Gen}T=T^{\bot}$, $T$ is called a large tilting sheaf.

{\bf{Theorem 6.15}}  Let $q\in\mathbb{Q}\cup\{\infty\}$, $G_{q}\oplus(\oplus_{S_{q}} S_{q}[\infty])$ is a large tilting sheaf, where $S_{q}$ runs through all quasi-simple sheaves of slope $q$.

{\bf{Proof}} Let $T=G_{q}\oplus(\oplus_{S_{q}} S_{q}[\infty])$, it need to prove that $\mbox{Gen}T=T^{\bot}$. By Lemma 3.2, $\mbox{Ext}^{1}(M, N)=0$, where $M$, $N$ lie in $\{S_{q}[\infty], G_{q}\ |\ S_{q}$ runs through all quasi-simple sheaves of slope $q\}$. Thus $\mbox{Ext}^{1}(T, \mbox{Add}T)=0$. Since $\mbox{Qcoh}(\mathbb{X})$ is hereditary, $\mbox{Gen}T\subset T^{\bot}$. Given any $F\in T^{\bot}$, by Proposition 6.1, denote $F=F_{1}\oplus F_{2}$, where $F_{1}\in\mathcal{C}_{q}$, $F_{2}\in\mathcal{Q}_{q}$. $\mbox{Ext}^{1}(T, F)=0$, so for each $S_{q}$ there has $\mbox{Ext}^{1}(S_{q}, F)=0$, i.e. $F$ is $q$-divisible. So $F_{1}$ lies in $\omega_{q}$, by Theorem 6.4, $F_{1}\in\mbox{Gen}T$. By Theorem 6.10
$F_{2}\in\mbox{Gen}T$, $F\in\mbox{Gen}T$. In summary, $\mbox{Gen}T=T^{\bot}$, $T$ is a large tilting sheaf.$\hfill\square$

{\bf Acknowledgements} {\rm The authors would like to thank Professor
Helmut Lenzing for his lectures about weighted projective lines in
Xiamen University in 2011. }

\bibliographystyle{amsplain}

\end{document}